\documentclass{article}
\usepackage{arxiv}
\usepackage[utf8]{inputenc} 
\usepackage[T1]{fontenc}    
\usepackage{hyperref}       
\usepackage{url}            
\usepackage{booktabs}       
\usepackage{amsfonts}       
\usepackage{nicefrac}       
\usepackage{microtype}      
\usepackage{lipsum}		
\usepackage{graphicx}
\usepackage{mathrsfs}
\usepackage{doi}
\usepackage{amsmath}
\usepackage{algorithm}
\usepackage{algorithmicx}
\usepackage[noend]{algpseudocode}
\usepackage{caption,subcaption}
\allowdisplaybreaks
\algnewcommand{\LeftComment}[1]{\Statex \(\triangleright\) #1}
\usepackage{stmaryrd}
\usepackage{color,soul}
\usepackage{wrapfig}
\usepackage{cite}
\hypersetup{
    colorlinks=true,
    linkcolor=blue,
    filecolor=magenta,      
    urlcolor=cyan,
}
\usepackage{courier}

\usepackage{listings}
\usepackage{xcolor}
\definecolor{codegreen}{rgb}{0,0.6,0}
\definecolor{codegray}{rgb}{0.5,0.5,0.5}
\definecolor{codepurple}{rgb}{0.58,0,0.82}
\definecolor{backcolour}{rgb}{0.95,0.95,0.92}

\lstdefinestyle{mystyle}{
  backgroundcolor=\color{backcolour},   commentstyle=\color{codegreen},
  keywordstyle=\color{magenta},
  numberstyle=\tiny\color{codegray},
  stringstyle=\color{codepurple},
  basicstyle=\ttfamily\footnotesize,
  breakatwhitespace=false,         
  breaklines=true,                 
  captionpos=b,                    
  keepspaces=true,                 
  numbers=left,                    
  numbersep=5pt,                  
  showspaces=false,                
  showstringspaces=false,
  showtabs=false,                  
  tabsize=2
}

\title{Augmented Lagrangian for hanging nodes in hexahedral meshes}

\author{
  {
  Saumik Dana}\\
	University of Southern California\\
	Los Angeles, CA 90007 \\
	\texttt{sdana@usc.edu} \\
}

\begin{document}
\maketitle
\begin{abstract}
The surge of activity in the resolution of fine scale features in the field of earth sciences over the past decade necessitates the development of robust yet simple algorithms that can tackle the various drawbacks of in silico models developed hitherto. One such drawback is that of the restrictive computational cost of finite element method in rendering resolutions to the fine scale features while at the same time keeping the domain being modeled sufficiently large. We propose the use of the augmented lagrangian commonly used in the treatment of hanging nodes in contact mechanics in tackling the drawback. An interface is introduced in a general hexahedral finite element mesh across which an aggressive coarsening of the finite elements is possible. The method is based upon minimizing an augmented potential energy which factors in the constraint that exists at the hanging nodes on that interface. This allows for a significant reduction in the number of finite elements comprising the mesh with concomitant reduction in the computational expense. 
\end{abstract}
\section{Introduction}
\begin{figure}[htb!]
\centering
\includegraphics[scale=0.45]{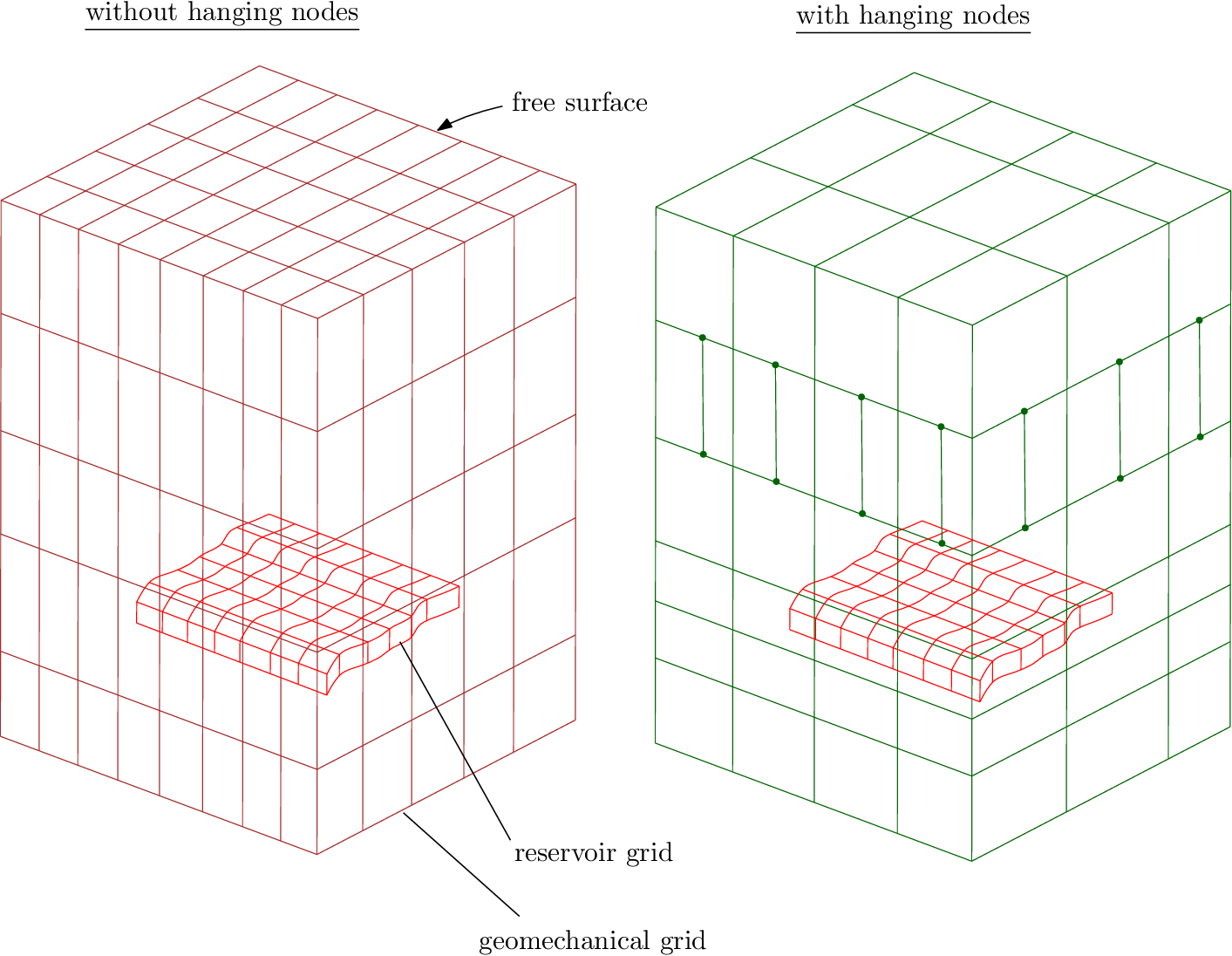}
\caption{The method in \cite{dana-2018} allows coarse grid for geomechanics coupled with fine grid for flow as shown on the left. The presence of hanging nodes in geomechanics grid as shown on the right allows the capability of capturing fine scale geomechanical features. The hanging nodes are represented by black dots to the right.}
\label{hanginggeo}
\end{figure}
The quantum of work devoted to modeling of fine scale features in the subsurface in the recent decade has spawned a need for simple yet powerful algorithms to simulate the same in silico with low computational cost. The main barrier to these simulations lies in the restrictively fine mesh that needs to be invoked to resolve the finer features of the corresponding physics while at the same time keeping the domain under consideration sufficiently large. The most logical approach to this problem is to allow for a fine mesh to exist in the regions which need a fine mesh and a coarse mesh to exist in regions which do not need a fine mesh. The authors previously developed a method to simulate subsurface flow on a fine mesh and subsurface mechanics on a coarse mesh while allowing for the coupling between the physics of flow and mechanics via a staggered solution algorithm~\cite{dana-2018}. The aforementioned work though is restrictive in the sense that the mesh for the mechanics domain needs to be uniformly coarser than the mesh for the flow domain as shown in Figure \ref{hanginggeo}. This makes the algorithm infeasible for problems involving fine scale features for the mechanics. With that in mind, we propose an addendum to the algorithm of \cite{dana-2018} by invoking the concept of hanging nodes in finite elements~\cite{Carlo-1978,M-1978,J.O-1985,Jua-1985,P-1985,H-1989,Panayioti-1992,P-1993,T-2000,nag-2001,becker-2003,puso-2003,Michae-2004,wriggersbook,P-2008} and the augmented lagrangian method~\cite{simo1992augmented,glowinski1989augmented,adeli1994augmented,conn1991globally} for treatment of hanging nodes. A depiction of geomechanics mesh with hanging nodes is given in Figure \ref{hanginggeo}. 
The problem is looked upon as minimization of a functional $\mathscr{C}$ with a constraint $\mathbf{g}=\mathbf{0}$ which dictates the geometry of the interface of the hanging nodes. The penalty formulation is
\begin{equation*}
\label{penalty}
\left.
\begin{array}{l}
Minimize\,\,\tilde{\mathscr{C}}\equiv \mathscr{C}+\frac{\epsilon}{2}\mathbf{g}\cdot \mathbf{g}\\
Subject\,\,to\,\,\mathbf{g}=\mathbf{0}
\end{array}
\right\|\qquad \mathrm{Penalty\,\,formulation}
\end{equation*}
where $\epsilon$ is a penalty parameter. A large enough $\epsilon$ lends to more accuracy while at the same time leading to highly ill-conditioned stiffness matrix in the eventual system of equations obtained at the discrete level. As a result, the choice of $\epsilon$ is a compromise between solution accuracy and solution stability. 
The lagrangian formulation is
\begin{equation*}
\label{lagrange}
\left.\begin{array}{l}
Minimize\,\,\tilde{\mathscr{C}}\equiv \mathscr{C}+\boldsymbol{\lambda}\cdot \mathbf{g}\\
Subject\,\,to\,\,\mathbf{g}=\mathbf{0}
\end{array}\right\|\qquad \mathrm{Lagrangian\,\,formulation}
\end{equation*}
where $\boldsymbol{\lambda}$ is the force conjugate to the constraint and is refered to as the lagrange multiplier. Although this method allows for the exact satisfaction of the constraint, the increase in number of degrees of freedom of the original system by the number of lagrange multipliers makes the augmentation computationally expensive. 
The perturbed Lagrangian formulation is
\begin{equation*}
\label{perturbed}
\left.\begin{array}{l}
Minimize\,\,\tilde{\mathscr{C}}\equiv \mathscr{C}+\boldsymbol{\lambda}\cdot \mathbf{g}-\frac{1}{2\epsilon}\boldsymbol{\lambda} \cdot \boldsymbol{\lambda}\\
Subject\,\,to\,\,\mathbf{g}-\frac{\boldsymbol{\lambda}}{\epsilon}=\mathbf{0}
\end{array}\right\|\qquad \mathrm{Perturbed\,\,lagrangian\,\,formulation}
\end{equation*}
This allows for the lagrange multiplier to be posed in terms of the constraint thus negating the need to solve for the multiplier as an additional degree of freedom. This method, though, suffers from the same problem that the original penalty method suffers from, i.e. a careful compromise between accuracy and stability must be made in the choice of the penalty parameter. 
The augmented Lagrangian formulation is 
\begin{equation*}
\label{augmented}
\left.\begin{array}{l}
Minimize\,\,\tilde{\mathscr{C}}\equiv \mathscr{C}+\boldsymbol{\lambda}^k\cdot \mathbf{g}+\frac{\epsilon}{2}\mathbf{g} \cdot \mathbf{g}\\
Subject\,\,to\,\,\boldsymbol{\lambda}^{k+1}-\boldsymbol{\lambda}^{k}=\epsilon \mathbf{g}
\end{array}\right\|\qquad \mathrm{Augmented\,\,lagrangian\,\,formulation}
\end{equation*}
where $\boldsymbol{\lambda}^{k}$ is the lagrange multiplier evaluated at the $k^{th}$ iteration. As is evident from the formulation, the lagrange multiplier is evaluated iteratively till it reaches an asymptotic value. The lagrange multiplier, is not an additional degree of freedom, and hence the system size does not increase as compared to the original minimization problem. The biggest advantage of this method is that the solution stability is not a function of the penalty parameter, and furthermore the lagrange multiplier iterative process reaches the true asymptotic value regardless of the value of the penalty parameter.
\section{Formulation}
\begin{figure}[htb!]
\centering
\includegraphics[scale=0.75]{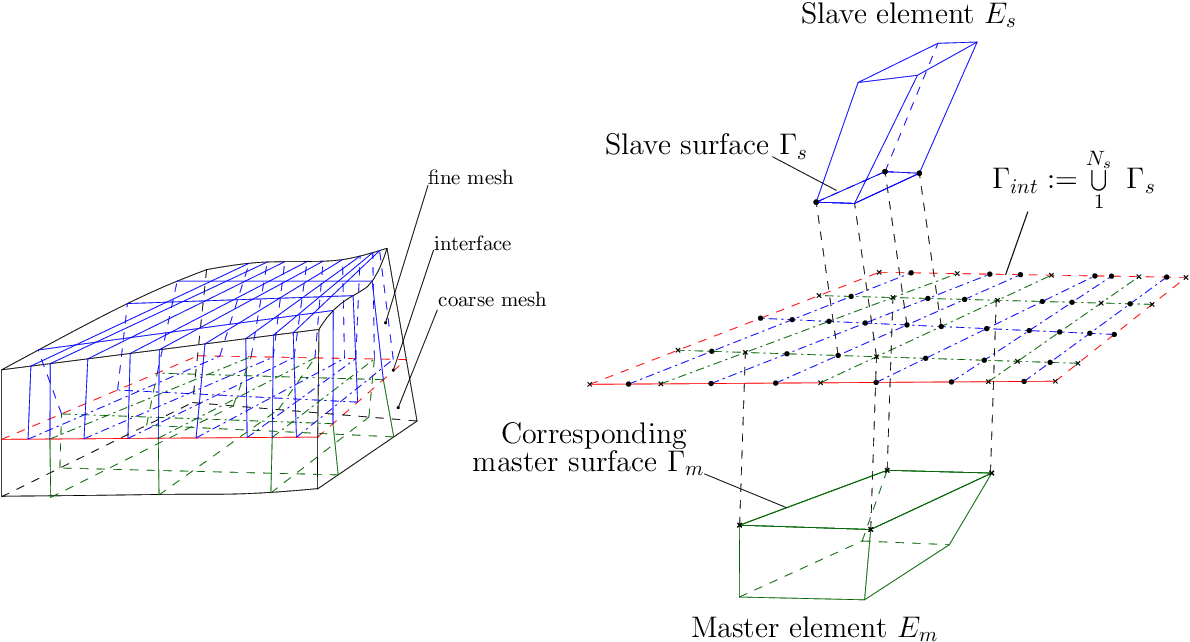}
\caption{There is an interface in the mesh across which an aggressive refinement is possible thus allowing for fine elements on one side of the interface and coarser elements on the other side of the interface}
\label{depict1}
\end{figure}
As shown in Figure \ref{depict1}, the presence of hanging nodes essentially means that there is an interface in the mesh across which an aggressive refinement is possible thus allowing for fine elements on one side of the interface and coarser elements on the other side of the interface. The fine and coarse elements are refered to as `slave element' and `master element' respectively while the faces of the slave and master elements making up the interface are refered to as `slave surface' and `master surface' respectively. Let $\mathbf{u}_s$ and $\mathbf{u}_m$ represent the displacement fields evaluated at $\Gamma_s$ and $\Gamma_m$ respectively. Then the problem statement is
\begin{equation}
\label{functional}
\left.\begin{array}{c}
Minimize\,\,\tilde{\mathscr{C}}\equiv \mathscr{C} + \sum\limits_{N_s}\int\limits_{\Gamma_s}\boldsymbol{\lambda}\cdot \mathbf{g}\,dA + \sum\limits_{N_s}\frac{1}{2}\int\limits_{\Gamma_s}\epsilon\mathbf{g}\cdot \mathbf{g}\, dA\\
Subject\,\,to\,\,\mathbf{g}\equiv \mathbf{u}_s-\mathbf{u}_m=\mathbf{0}\qquad \forall\,\,\Gamma_s
\end{array}\right.
\end{equation}
where $\mathscr{C}$ is the strain energy in the absence of hanging nodes, $\mathbf{g}$ is the refered to as the penetration function, 
\begin{align*}
\sum\limits_{N_s}\int\limits_{\Gamma_s}\boldsymbol{\lambda}\cdot \mathbf{g}\,dA
\end{align*}
is the lagrange multiplier term with $\boldsymbol{\lambda}$ being the lagrange multiplier and 
\begin{align*}
\sum\limits_{N_s}\frac{1}{2}\int\limits_{\Gamma_s}\epsilon\mathbf{g}\cdot \mathbf{g}\, dA
\end{align*}
is the penalty term with $\epsilon$ being the penalty parameter. Let $\mathbf{t}^s$ and $\mathbf{t}^m$ be force conjugates to the constraint $\mathbf{g}=\mathbf{0}$ at $\Gamma_s$ and $\Gamma_m$ respectively. Then 
\begin{align*}
\boldsymbol{\lambda}\equiv\frac{1}{2}(\mathbf{t}^s+\mathbf{t}^m)
\end{align*}
is the force conjugate to the constraint $\mathbf{g}=\mathbf{0}$ introduced in a mean sense. 

For the sake of clarity, we rewrite $\tilde{\mathscr{C}}$ as
\begin{align}
\label{fagain}
\tilde{\mathscr{C}}\equiv \mathscr{C} + \sum\limits_{N_s}\int\limits_{\Gamma_s}\frac{1}{2}(\mathbf{t}^s+\mathbf{t}^m)\cdot (\mathbf{u}_s-\mathbf{u}_m)\,dA 
+\sum\limits_{N_s}\frac{1}{2}\int\limits_{\Gamma_s}\epsilon(\mathbf{u}_s-\mathbf{u}_m)\cdot (\mathbf{u}_s-\mathbf{u}_m)\, dA
\end{align}
Minimization of \eqref{fagain} would imply equating the first variation to zero as follows
\begin{equation}
\label{mini}
\left.\begin{array}{l}
\delta \tilde{\mathscr{C}}\equiv 0 =\delta \mathscr{C}+ \mathcal{C}
\end{array}\right.
\end{equation}
where $\mathcal{C}$ is given by
\begin{equation}
\label{mini1}
\left.\begin{array}{l}
\mathcal{C}:= \sum\limits_{N_s}\int\limits_{\Gamma_s}\frac{1}{2}(\delta \mathbf{t}^s+\delta \mathbf{t}^m)\cdot (\mathbf{u}_s-\mathbf{u}_m)\,dA\\
\quad + \sum\limits_{N_s}\int\limits_{\Gamma_s}\frac{1}{2}(\mathbf{t}^s+\mathbf{t}^m)\cdot (\delta \mathbf{u}_s-\delta \mathbf{u}_m)\,dA\\
\quad +\sum\limits_{N_s}\int\limits_{\Gamma_s}\epsilon (\mathbf{u}_s-\mathbf{u}_m) \cdot (\delta \mathbf{u}_s-\delta \mathbf{u}_m)\,dA
\end{array}\right.
\end{equation}

\begin{figure}[htb!]
\centering
\includegraphics[trim={0 0 0 5.5cm},clip,scale=0.5]{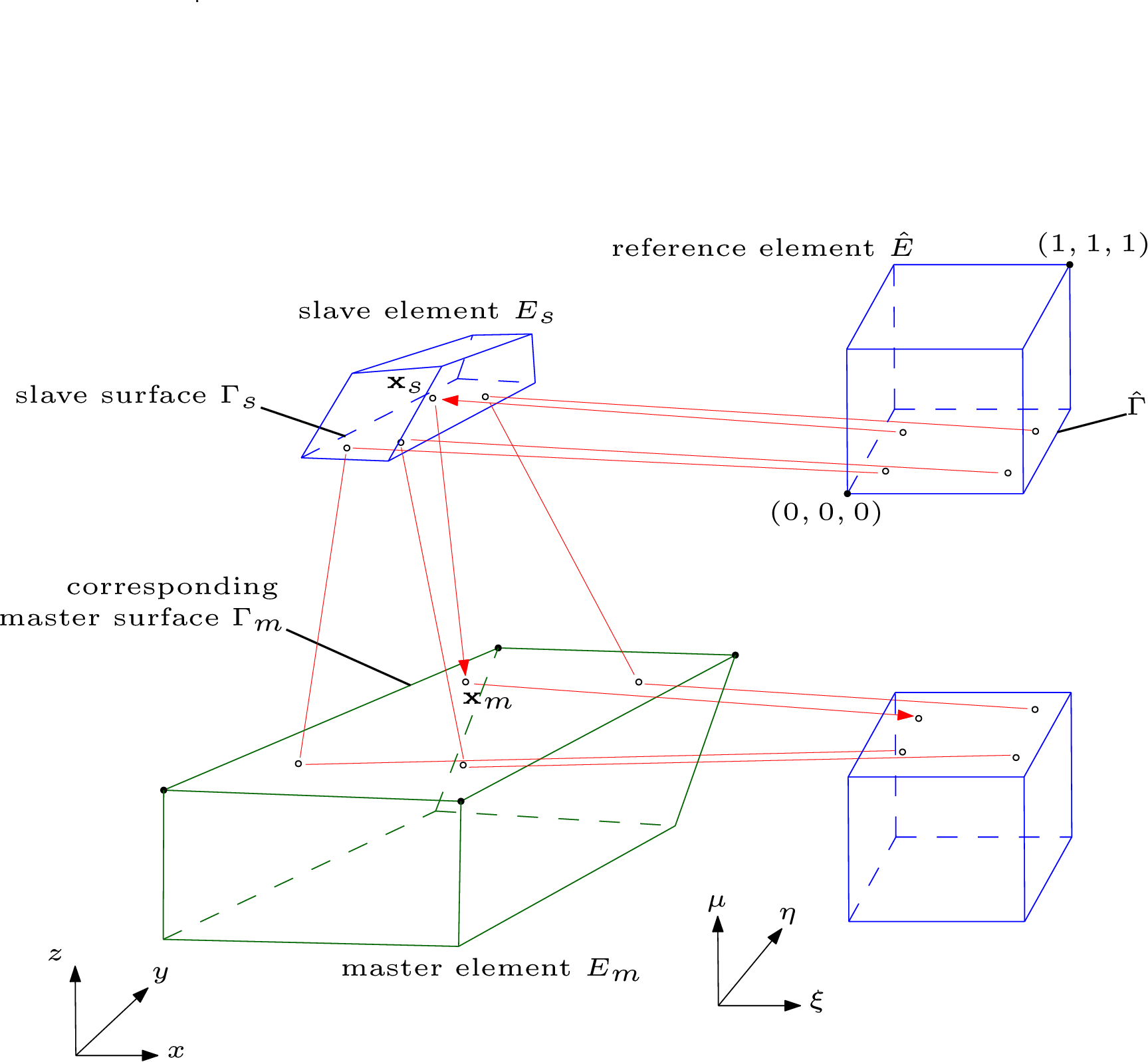}
\caption{Corresponding to each of the four gauss points on $\hat{\Gamma}$, there is an actual physical point on $\Gamma_s$. That point has an orthogonal projection onto $\Gamma_m$. That orthogonally projected point has a corresponding point on $\hat{E}$.}
\label{damn}
\end{figure}
The contribution to $\mathcal{C}$ over every $\Gamma_s$ is evaluated as a sum of the integrand $\hat{\mathcal{C}}$ evaluated at each of the four gauss points $g\in \mathcal{G}$ shown in Figure \ref{damn} multiplied by the determinant $J_{\Gamma_s}$ of the jacobian of the mapping $\hat{\Gamma}\mapsto \Gamma_s$ as follows
\begin{equation}
\label{mini2}
\left.\begin{array}{l}
\mathcal{C}:= \sum\limits_{\Gamma_s \in \Gamma_{int}}\sum\limits_{g\in \mathcal{G}}\frac{1}{2}(\delta \mathbf{t}^s+\delta \mathbf{t}^m)\cdot (\mathbf{u}_s-\mathbf{u}_m)\,J_{\Gamma_s}\\
\quad + \sum\limits_{\Gamma_s \in \Gamma_{int}}\sum\limits_{g\in \mathcal{G}}\frac{1}{2}(\mathbf{t}^s+\mathbf{t}^m)\cdot (\delta \mathbf{u}_s-\delta \mathbf{u}_m)\,J_{\Gamma_s}\\
\quad +\sum\limits_{\Gamma_s \in \Gamma_{int}}\sum\limits_{g\in \mathcal{G}}\epsilon (\mathbf{u}_s-\mathbf{u}_m) \cdot (\delta \mathbf{u}_s-\delta \mathbf{u}_m)\,J_{\Gamma_s}
\end{array}\right.
\end{equation}

\section{System of equations}
As shown in Figure \ref{damn}, corresponding to each gauss point $(\xi_s,\eta_s,-1)$ on $\hat{\Gamma}$, there is an actual physical point $\mathbf{x}_s$ on $\Gamma_s$ given by
\begin{equation}
\label{coord}
\mathbf{x}_s = \sum\limits_{i=1}^8 N_i\lvert_{(\xi_s,\eta_s,-1)}\mathbf{X}_{s}^i\equiv \mathcal{N}_{E_s}\mathbf{X}_{s}
\end{equation}
where $\mathbf{X}_{s}^i,\,i=1,..,8$ are coordinates of nodes of $E_s$ and $N_i(\xi,\eta),\,i=1,..,8$ represent the shape functions.
Let $\mathbf{x}_m$ be the orthogonal projection of $\mathbf{x}_s$ onto the corresponding master surface with corresponding location $\chi\equiv (\xi_m,\eta_m,1)$ on $\hat{E}$ such that
\begin{equation}
\label{coord1}
\mathbf{x}_m = \sum\limits_{i=1}^8 N_i\lvert_{\chi}\mathbf{X}_{m}^i\equiv \mathcal{N}_{E_m}\mathbf{X}_{m}
\end{equation}
where $\mathbf{X}_{m}^i,\,i=1,..,8$ be the coordinates of nodes of $E_m$. 
We know $\mathbf{x}_s$ but need to evaluate $\mathbf{x}_m$.
\subsection{Evaluating $\mathbf{x}_m$ given $\mathbf{x}_s$}
The orthogonality condition is satisfied by 
\begin{equation}
\label{thirdf}
\left.\begin{array}{c}
\mathbf{e}_1\cdot (\mathbf{x}_s-\mathbf{x}_m)=0,\qquad \mathbf{e}_2\cdot (\mathbf{x}_s-\mathbf{x}_m)=0, \qquad \mathbf{e}_3\cdot (\mathbf{x}_s-\mathbf{x}_m)=0
\end{array}\right.
\end{equation}
where the components $\mathbf{e}_1$, $\mathbf{e}_2$ and $\mathbf{e}_3$ of the tangent at $\mathbf{x}_m$ with respect to the local axis of master surface are computed as
\begin{equation}
\label{secondf}
\left.\begin{array}{c}
\mathbf{e}_1 = \sum\limits_{j=1}^8\frac{\partial N_j}{\partial \xi}\lvert_{\chi}\mathbf{X}_{m_j},\qquad 
\mathbf{e}_2 = \sum\limits_{j=1}^8\frac{\partial N_j}{\partial \eta}\lvert_{\chi}\mathbf{X}_{m_j}\qquad 
\mathbf{e}_3 = \sum\limits_{j=1}^8\frac{\partial N_j}{\partial \mu}\lvert_{\chi}\mathbf{X}_{m_j}
\end{array}\right.
\end{equation}
Substituting \eqref{secondf}, \eqref{coord} and \eqref{coord1} in \eqref{thirdf}, we get
\begin{equation}
\label{fourthf}
\left.\begin{array}{c}
\sum\limits_{j=1}^8 \frac{\partial N_j}{\partial \xi}\lvert_{\chi}\mathbf{X}_{m_j}\cdot\bigg(\mathbf{x}_s -\sum\limits_{i=1}^8 N_i\lvert_{\chi} \mathbf{X}_{m_i}\bigg)\equiv f_1\vert_{\chi}=0\\
\sum\limits_{j=1}^8 \frac{\partial N_j}{\partial \eta}\lvert_{\chi}\mathbf{X}_{m_j}\cdot\bigg(\mathbf{x}_s -\sum\limits_{i=1}^8 N_i\lvert_{\chi} \cdot \mathbf{X}_{m_i}\bigg)\equiv f_2\vert_{\chi}=0\\
\sum\limits_{j=1}^8 \frac{\partial N_j}{\partial \mu}\lvert_{\chi}\mathbf{X}_{m_j}\cdot\bigg(\mathbf{x}_s -\sum\limits_{i=1}^8 N_i\lvert_{\chi} \cdot \mathbf{X}_{m_i}\bigg)\equiv f_3\vert_{\chi}=0
\end{array}\right.
\end{equation}
The solution to \eqref{fourthf} is obtained iteratively for the $(k+1)^{th}$ iteration as 
\begin{equation*}
\chi^{k+1}
=\chi^{k} 
-\begin{bmatrix}
\frac{\partial f_1}{\partial \xi}\vert_{\chi^{k}} & \frac{\partial f_1}{\partial \eta}\vert_{\chi^{k}} & \frac{\partial f_1}{\partial \mu}\vert_{\chi^{k}}\\
\frac{\partial f_2}{\partial \xi}\vert_{\chi^{k}} & \frac{\partial f_2}{\partial \eta}\vert_{\chi^{k}} & \frac{\partial f_2}{\partial \mu}\vert_{\chi^{k}}\\
\frac{\partial f_3}{\partial \xi}\vert_{\chi^{k}} & \frac{\partial f_3}{\partial \eta}\vert_{\chi^{k}} & \frac{\partial f_3}{\partial \mu}\vert_{\chi^{k}}
\end{bmatrix}^{-1}
\left\{
\begin{array}{c}
f_1\vert_{\chi^k}\\
f_2\vert_{\chi^k}\\
f_3\vert_{\chi^k}
\end{array}
\right\}
\end{equation*}
with initial guess 
\begin{align*}
\chi^{0}=\left\{\begin{array}{c}0\\0\\1\end{array}\right\}
\end{align*}
The stopping criterion is
\begin{equation*}
\Vert \chi^{k+1}-\chi^{k}\Vert< TOL* \Vert\chi^{k}\Vert
\end{equation*}
where $TOL$ is a pre-specified tolerance. Once this criterion is satisfied, we set $\chi=\chi^{k+1}$ and then obtain $\mathbf{x}_m$ using \eqref{coord1}.

\subsection{Evaluating $\mathbf{u}_s$, $\delta \mathbf{u}_s$, $\mathbf{u}_m$ and $\delta \mathbf{u}_m$; $\mathbf{t}^s$, $\delta \mathbf{t}^s$, $\mathbf{t}^m$ and $\delta \mathbf{t}^m$}

Let $\mathbf{U}$ represent the vector of nodal displacement degrees of freedom, and let $\mathbf{U}\vert_{E}$ represent the restriction of $\mathbf{U}$ to any element $E$. Then we have
\begin{equation}
\label{pehla}
\left.\begin{array}{c}
\mathbf{u}_s = \sum\limits_{i=1}^8 N_i\lvert_{(\xi_s,\eta_s,-1)}\mathbf{U}\vert_{E_s}^i\equiv \mathcal{N}_{E_s}\mathbf{U}\vert_{E_s},\qquad 
\mathbf{u}_m = \sum\limits_{i=1}^8 N_i\lvert_{(\xi_m,\eta_m,1)}\mathbf{U}\vert_{E_m}^i\equiv \mathcal{N}_{E_m}\mathbf{U}\vert_{E_m}
\end{array}\right.
\end{equation}
The force conjugate to the constraint evaluated at $\mathbf{x}_s$ is given by
\begin{align}
\nonumber
\mathbf{t}^s &= 
\begin{bmatrix}
\sigma_1 & \sigma_4 & \sigma_6\\
\sigma_4 & \sigma_2 & \sigma_5\\
\sigma_6 & \sigma_5 & \sigma_3
\end{bmatrix}\bigg\vert_{\mathbf{x}_s}
\left\{\begin{array}{c}
n_1\\
n_2\\
n_3
\end{array}
\right\}\bigg\vert_{\mathbf{x}_s}
\equiv 
\begin{bmatrix}
n_1 & 0 & 0 & n_2 & 0 & n_3\\
0 & n_2 & 0 & n_1 & n_3 & 0\\
0 & 0 & n_3 & 0 & n_2 & n_1
\end{bmatrix}\bigg\vert_{\mathbf{x}_s}
\left\{\begin{array}{c}
\sigma_1\equiv \sigma_{xx}\\
\sigma_2\equiv \sigma_{yy}\\
\sigma_3\equiv \sigma_{zz}\\
\sigma_4\equiv \sigma_{xy}\\
\sigma_5\equiv \sigma_{yz}\\
\sigma_6\equiv \sigma_{xz}
\end{array}
\right\}\bigg\vert_{\mathbf{x}_s}\\
\label{doosra}
&\equiv 
\overbrace{\begin{bmatrix}
n_1 & 0 & 0 & n_2 & 0 & n_3\\
0 & n_2 & 0 & n_1 & n_3 & 0\\
0 & 0 & n_3 & 0 & n_2 & n_1
\end{bmatrix}\bigg\vert_{\mathbf{x}_s}\mathcal{D}\mathbf{B}\vert_{\mathbf{x}_s}}^{\mathcal{F}_{E_s}}\mathbf{U}\vert_{E_s}
\end{align}
where $\left\{\begin{array}{c}
n_1\\
n_2\\
n_3
\end{array}
\right\}\bigg\vert_{\mathbf{x}_s}$ is the normal to $\Gamma_s$ evaluated at $\mathbf{x}_s$, $\mathcal{D}$ is the $6\times 6$ constitutive matrix and $\mathbf{B}\vert_{\mathbf{x}_s}$ is the $6\times 24$ strain displacement interpolation matrix evaluated at $\mathbf{x}_s$. 

Similarly, the force conjugate to the constraint evaluated at $\mathbf{x}_m$ is given by 
\begin{align}
\nonumber
\mathbf{t}^m &= 
\begin{bmatrix}
\sigma_1 & \sigma_4 & \sigma_6\\
\sigma_4 & \sigma_2 & \sigma_5\\
\sigma_6 & \sigma_5 & \sigma_3
\end{bmatrix}\bigg\vert_{\mathbf{x}_m}
\left\{\begin{array}{c}
n_1\\
n_2\\
n_3
\end{array}
\right\}\bigg\vert_{\mathbf{x}_m}
\equiv 
\begin{bmatrix}
n_1 & 0 & 0 & n_2 & 0 & n_3\\
0 & n_2 & 0 & n_1 & n_3 & 0\\
0 & 0 & n_3 & 0 & n_2 & n_1
\end{bmatrix}\bigg\vert_{\mathbf{x}_m}
\left\{\begin{array}{c}
\sigma_1\equiv \sigma_{xx}\\
\sigma_2\equiv \sigma_{yy}\\
\sigma_3\equiv \sigma_{zz}\\
\sigma_4\equiv \sigma_{xy}\\
\sigma_5\equiv \sigma_{yz}\\
\sigma_6\equiv \sigma_{xz}
\end{array}
\right\}\bigg\vert_{\mathbf{x}_m}\\
\label{teesra}
&\equiv 
\overbrace{\begin{bmatrix}
n_1 & 0 & 0 & n_2 & 0 & n_3\\
0 & n_2 & 0 & n_1 & n_3 & 0\\
0 & 0 & n_3 & 0 & n_2 & n_1
\end{bmatrix}\bigg\vert_{\mathbf{x}_m}\mathcal{D}\mathbf{B}\vert_{\mathbf{x}_m}}^{\mathcal{F}_{E_m}}\mathbf{U}\vert_{E_m}
\end{align}
where $\left\{\begin{array}{c}
n_1\\
n_2\\
n_3
\end{array}
\right\}\bigg\vert_{\mathbf{x}_m}$ is the normal to $\Gamma_m$ evaluated at $\mathbf{x}_m$ and $\mathbf{B}\vert_{\mathbf{x}_m}$ is the $6\times 24$ strain displacement interpolation matrix evaluated at $\mathbf{x}_m$. 

The normals $\left\{\begin{array}{c}
n_1\\
n_2\\
n_3
\end{array}
\right\}\bigg\vert_{\mathbf{x}_s}$ and $\left\{\begin{array}{c}
n_1\\
n_2\\
n_3
\end{array}
\right\}\bigg\vert_{\mathbf{x}_m}$ are obtained as follows
\begin{align*}
\left\{\begin{array}{c}
n_1\\
n_2\\
n_3
\end{array}
\right\}\bigg\vert_{\mathbf{x}_s} = \frac{\nabla S_s}{\lVert \nabla S_s\rVert}\bigg\vert_{\mathbf{x}_s},\qquad 
\left\{\begin{array}{c}
n_1\\
n_2\\
n_3
\end{array}
\right\}\bigg\vert_{\mathbf{x}_m} = \frac{\nabla S_m}{\lVert \nabla S_m\rVert}\bigg\vert_{\mathbf{x}_m}
\end{align*}
where $S_s$ and $S_m$ are equations of the slave and master surfaces respectively. The procedure to obtain equations of faces of the elements in given in Appendix~\ref{surface}.

\subsection{Evaluating the surface integral}

Let $\mathcal{E}_s$ be the collection of all slave elements. In lieu of Equations \eqref{pehla} - \eqref{teesra}, the surface integral \eqref{mini2} is evaluated as
\begin{align}
\nonumber
\mathcal{C}=\sum\limits_{\mathscr{E}_s\in \mathcal{E}_s}^{}\bigg[\sum\limits_{N=1}^4\bigg[&\frac{1}{2}(\mathcal{F}_{\mathscr{E}_s}\delta \mathbf{U}\vert_{\mathscr{E}_s}+\mathcal{F}_{\mathscr{E}_m}\delta \mathbf{U}\vert_{\mathscr{E}_m})\cdot (\mathcal{N}_{\mathscr{E}_s}\mathbf{U}\vert_{\mathscr{E}_s}
-\mathcal{N}_{\mathscr{E}_m}\mathbf{U}\vert_{\mathscr{E}_m})\\
\nonumber
 + &\frac{1}{2}(\mathcal{F}_{\mathscr{E}_s} \mathbf{U}\vert_{\mathscr{E}_s}+\mathcal{F}_{\mathscr{E}_m} \mathbf{U}\vert_{\mathscr{E}_m})\cdot (\mathcal{N}_{\mathscr{E}_s}\delta\mathbf{U}\vert_{\mathscr{E}_s}
-\mathcal{N}_{\mathscr{E}_m}\delta\mathbf{U}\vert_{\mathscr{E}_m})\\
\label{gandu}
+ &\epsilon (\mathcal{N}_{\mathscr{E}_s}\mathbf{U}\vert_{\mathscr{E}_s}
-\mathcal{N}_{\mathscr{E}_m}\mathbf{U}\vert_{\mathscr{E}_m}) \cdot  (\mathcal{N}_{\mathscr{E}_s}\delta \mathbf{U}\vert_{\mathscr{E}_s}
-\mathcal{N}_{\mathscr{E}_m}\delta \mathbf{U}\vert_{\mathscr{E}_m}) \bigg]\,det J_{\mathscr{E}_s}\bigg]
\end{align}
Which can also be written as
\begin{align}
\nonumber
\mathcal{C}&=\delta \mathbf{U}_s^T \bigg[\overbrace{\sum\limits_{\mathscr{E}_s\in \mathcal{E}_s}^{}\bigg[\overbrace{\sum\limits_{N=1}^4\bigg[\frac{1}{2}\mathcal{F}_{\mathscr{E}_s}^T\mathcal{N}_{\mathscr{E}_s}+\frac{1}{2}\mathcal{N}_{\mathscr{E}_s}^T\mathcal{F}_{\mathscr{E}_s}+\epsilon \mathcal{N}_{\mathscr{E}_s}^T \mathcal{N}_{\mathscr{E}_s}\bigg]det J_{\mathscr{E}_s}}^{\mathbf{K}_{ss\vert_{\mathscr{E}_s}}}\bigg]}^{\mathbf{K}_{ss}}\bigg]\mathbf{U}_s\\
\nonumber
&+ \delta \mathbf{U}_s^T\bigg[\overbrace{\sum\limits_{\mathscr{E}_s\in \mathcal{E}_s}^{}\bigg[\overbrace{\sum\limits_{N=1}^4\bigg[-\frac{1}{2}\mathcal{F}_{\mathscr{E}_s}^T\mathcal{N}_{\mathscr{E}_m}+\frac{1}{2}\mathcal{N}_{\mathscr{E}_s}^T\mathcal{F}_{\mathscr{E}_m}-\epsilon \mathcal{N}_{\mathscr{E}_s}^T \mathcal{N}_{\mathscr{E}_m}\bigg]det J_{\mathscr{E}_s}}^{\mathbf{K}_{sm\vert_{\mathscr{E}_s}}}\bigg]}^{\mathbf{K}_{sm}}\bigg]\mathbf{U}_m\\ 
\nonumber
&+\delta \mathbf{U}_m^T\bigg[\overbrace{\sum\limits_{\mathscr{E}_s\in \mathcal{E}_s}^{}\bigg[\overbrace{\sum\limits_{N=1}^4\bigg[\frac{1}{2}\mathcal{F}_{\mathscr{E}_m}^T\mathcal{N}_{\mathscr{E}_s}-\frac{1}{2}\mathcal{N}_{\mathscr{E}_m}^T\mathcal{F}_{\mathscr{E}_s}-\epsilon \mathcal{N}_{\mathscr{E}_m}^T \mathcal{N}_{\mathscr{E}_s}\bigg]det J_{\mathscr{E}_s}}^{\mathbf{K}_{ms\vert_{\mathscr{E}_s}}}\bigg]}^{\mathbf{K}_{ms}}\bigg]\mathbf{U}_s\\
\label{dabba}
&+ \delta \mathbf{U}_m^T\bigg[\overbrace{\sum\limits_{\mathscr{E}_s\in \mathcal{E}_s}^{}\bigg[\overbrace{\sum\limits_{N=1}^4\bigg[-\frac{1}{2}\mathcal{F}_{\mathscr{E}_m}^T\mathcal{N}_{\mathscr{E}_m}-\frac{1}{2}\mathcal{N}_{\mathscr{E}_m}^T\mathcal{F}_{\mathscr{E}_m}+\epsilon \mathcal{N}_{\mathscr{E}_m}^T \mathcal{N}_{\mathscr{E}_m}\bigg]det J_{\mathscr{E}_s}}^{\mathbf{K}_{mm\vert_{\mathscr{E}_s}}}\bigg]}^{\mathbf{K}_{mm}}\bigg]\mathbf{U}_m
\end{align}
where $\mathbf{U}_s$ and $\mathbf{U}_m$ are the collection of displacement degrees of freedom corresponding to nodes of slave elements and master elements respectively.
The system of equations 
is eventually written as
\begin{align}
\nonumber
\begin{bmatrix}
\mathbf{K}^d+\begin{bmatrix}
. & . & .\\
. & \mathbf{K}_{ss} & \mathbf{K}_{sm}\\
. & \mathbf{K}_{ms} & \mathbf{K}_{mm}
\end{bmatrix}\end{bmatrix}\left\{\begin{array}{c}
\mathbf{U}_r\\
\mathbf{U}_s\\
\mathbf{U}_m
\end{array}\right\}=\mathbf{P}
\end{align}
where $\mathbf{U}_r$ is the collection of displacement degrees of freedom corresponding to nodes of all elements which are neither slave elements nor master elements, and $\mathbf{K}_{ss}$, $\mathbf{K}_{sm}$, $\mathbf{K}_{ms}$ and $\mathbf{K}_{mm}$ are given in Equation \eqref{dabba}.
\section{Procedural framework}
The steps to be followed for the treatment of hanging nodes in hexahedral meshes are 
\begin{itemize}
\item Identify the elements sharing the interface
\item Identify the elements on the fine mesh side as slave elements and elements on the coarse mesh side as master elements
\item Identify the faces of the slave elements on the interface as slave surfaces and faces of the master elements on the interface as master surfaces
\item Use singular value decompositions~\cite{dana-2018} to obtain the equations of the slave and master surfaces
\item In the numerical integration module, map the slave and master surfaces to 2D reference elements
\item For every gauss point on the reference element which every slave surface has been mapped onto, identify the point on the slave surface. 
\item Use the equation of the slave surface to obtain the normal to the slave surface at that point.
\item Obtain the orthogonal projection of that point onto the master surface.
\item Use the equation of the master surface to obtain the normal to the master surface at that point.
\item Obtain the contributions to the submatrices from each slave element
\item Assemble the contributions to obtain the global submatrices
\end{itemize}

\appendix
\section{Obtaining equations of the element faces}\label{surface}
\begin{figure}[htb!]
\centering
\includegraphics[scale=0.9]{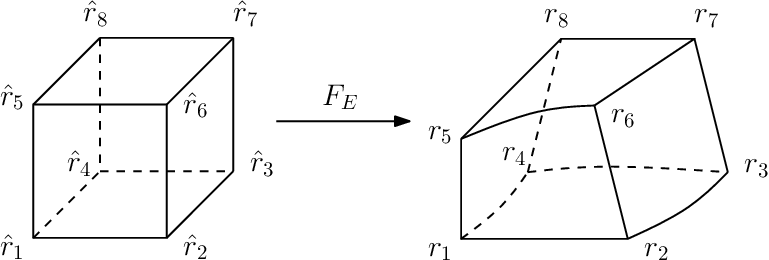}
\caption{Trilinear mapping $F_E:\hat{E}\rightarrow E$ for 8 noded distorted hexahedral elements. The faces of $E$ can be non-planar.}
\label{mapping}
\end{figure}
Let $\mathcal{T}_h$ be finite element partition of $\Omega\subset \mathbb{R}^3$ consisting of distorted hexahedral elements $E$ where $h=\max_{E\in\mathcal{T}_h}diam(E)$. 
Let $\mathbf{r}_i$, $i=1,..,8$ be the vertices of $E$.
Now consider a reference cube $\hat{E}$ with vertices $\hat{\mathbf{r}}_1=[0\,\,0\,\,0]^T$, $\hat{\mathbf{r}}_2=[1\,\,0\,\,0]^T$, $\hat{\mathbf{r}}_3=[1\,\,1\,\,0]^T$, $\hat{\mathbf{r}}_4=[0\,\,1\,\,0]^T$, $\hat{\mathbf{r}}_5=[0\,\,0\,\,1]^T$, $\hat{\mathbf{r}}_6=[1\,\,0\,\,1]^T$, $\hat{\mathbf{r}}_7=[1\,\,1\,\,1]^T$ and $\hat{\mathbf{r}}_8=[0\,\,1\,\,1]^T$ as shown in Figure \ref{mapping}.
Let $\hat{\mathbf{x}}=(\hat{x},\hat{y},\hat{z})\in \hat{E}$ and $\mathbf{x}=(x,y,z)\in E$. The function $F_E(\hat{\mathbf{x}}):\hat{E}\rightarrow E$ is
\begin{align*} 
\nonumber F_E(\hat{\mathbf{x}})=\mathbf{r}_1(1-\hat{x})(1-\hat{y})(1-\hat{z})+\mathbf{r}_2\hat{x}(1-\hat{y})(1-\hat{z})+\mathbf{r}_3\hat{x}\hat{y}(1-\hat{z})+\mathbf{r}_4(1-\hat{x})\hat{y}(1-\hat{z})&\\
+\mathbf{r}_5(1-\hat{x})(1-\hat{y})\hat{z}+\mathbf{r}_6\hat{x}(1-\hat{y})\hat{z}
+\mathbf{r}_7\hat{x}\hat{y}\hat{z}+
\mathbf{r}_8(1-\hat{x})\hat{y}\hat{z}&
\end{align*}
Denote Jacobian matrix by $DF_E$ and let $J_E=det(DF_E)$. Defining $\mathbf{r}_{ij}\equiv \mathbf{r}_i-\mathbf{r}_j$, we have 
\begin{align*}
DF_E(\mathbf{\hat{x}}) = \begin{bmatrix}
\mathbf{r}_{21}+(\mathbf{r}_{34}-\mathbf{r}_{21})\hat{y}
+(\mathbf{r}_{65}-\mathbf{r}_{21})\hat{z}
+((\mathbf{r}_{21}-\mathbf{r}_{34})
-(\mathbf{r}_{65}-\mathbf{r}_{78}))\hat{y}\hat{z};\\
\mathbf{r}_{41}+(\mathbf{r}_{34}-\mathbf{r}_{21})\hat{x}
+(\mathbf{r}_{85}-\mathbf{r}_{41})\hat{z}
+((\mathbf{r}_{21}-\mathbf{r}_{34})
-(\mathbf{r}_{65}-\mathbf{r}_{78}))\hat{x}\hat{z};\\
\mathbf{r}_{51}+(\mathbf{r}_{65}-\mathbf{r}_{21})\hat{x}
+(\mathbf{r}_{85}-\mathbf{r}_{41})\hat{y}
+((\mathbf{r}_{21}-\mathbf{r}_{34})
-(\mathbf{r}_{65}-\mathbf{r}_{78}))\hat{x}\hat{y}
\end{bmatrix}_{3\times 3}
\end{align*}
Denote inverse mapping by $F_E^{-1}$, its Jacobian matrix by $DF_E^{-1}$ and let $J_{F_E^{-1}}=det(DF_E^{-1})$
such that
\begin{align*}
DF_E^{-1}(\mathbf{x})=(DF_E)^{-1}(\hat{\mathbf{x}});\qquad  J_{F_E^{-1}}(\mathbf{x})=(J_E)^{-1}(\hat{\mathbf{x}})
\end{align*}
Let $\phi(\mathbf{x})$ be any function defined on $E$ and $\hat{\phi}(\hat{\mathbf{x}})$ be its corresponding definition on $\hat{E}$. Then we have
\begin{align}
\label{hat2}
\nabla \phi = (DF_E^{-1})^T(\mathbf{x})\,\hat{\nabla} \hat{\phi} = (DF_E)^{-T}(\mathbf{\hat{x}})\,\hat{\nabla} \hat{\phi}
\end{align}

\begin{figure}[htb!]
\centering
\includegraphics[scale=0.8]{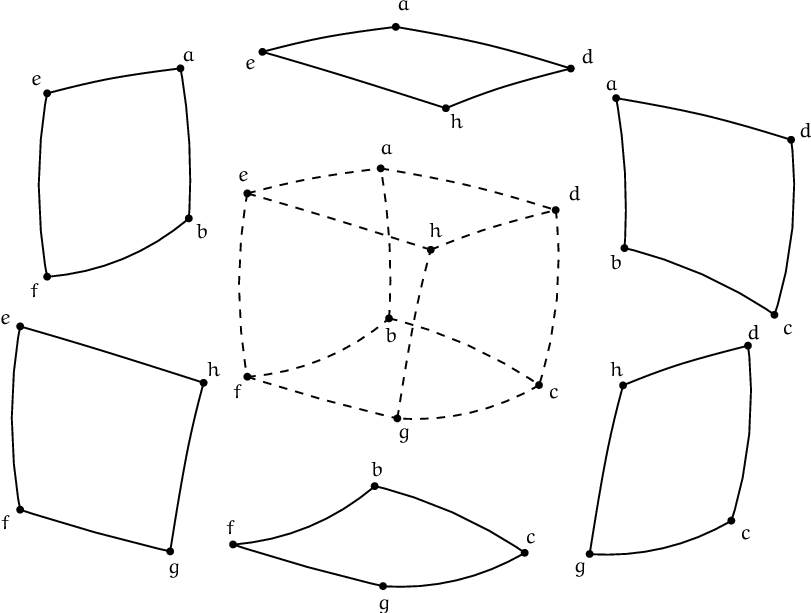}
\caption{A representation of hexahedral element $E\equiv abcdefgh$ with its six faces $aehd$, $abfe$, $ehgf$, $bcgf$, $cdhg$ and $adcb$. The coordinate information of the four vertices of each of the faces is used to obtain its equation.}
\label{gain1}
\end{figure}
\par
Let $\mathcal{S}(\mathbf{x})=0$, $\mathbf{x}\equiv(x,y,z)\in e$ be the equation of face $e$ of element $E$ with its vertices $\mathbf{v}_i\equiv(x_i,y_i,z_i)$, $i=1,2,3,4$. A representation of $E$ with its faces is provided in Figure \ref{gain1}.
Define $\mathcal{S}(\mathbf{x})$ by a trilinear as
\begin{align}
\label{expression}
\mathcal{S}(\mathbf{x})=\begin{bmatrix}xyz & xy & yz & xz & x & y & z & 1
\end{bmatrix} \mathbf{c}_{8 \times 1}
\end{align}
where $\mathbf{c}_{8\times 1}$ is the vector of coefficients to be determined. Since the equation $\mathcal{S}(\mathbf{x})=0$ is satisfied at each of the four vertices defining the face, we get the system of equations
\begin{align*}
\overbrace{\begin{bmatrix}
x_1y_1z_1 & x_1y_1 & y_1z_1 & x_1z_1 & x_1 & y_1 & z_1 & 1\\
x_2y_2z_2 & x_2y_2 & y_2z_2 & x_2z_2 & x_2 & y_2 & z_2 & 1\\
x_3y_3z_3 & x_3y_3 & y_3z_3 & x_3z_3 & x_3 & y_3 & z_3 & 1\\
x_4y_4z_4 & x_4y_4 & y_4z_4 & x_4z_4 & x_4 & y_4 & z_4 & 1
\end{bmatrix}}^{\mathbf{M}_{4\times 8}}\mathbf{c}_{8 \times 1}=
\begin{bmatrix}
0\\0\\0\\0
\end{bmatrix}_{4\times 1}
\end{align*}
for $\mathbf{c}$. The objective is to determine $\mathbf{c}\in Null(\mathbf{M})$. First, we get the SVD of $\mathbf{M}$ as
\begin{align}
\label{svd}
\mathbf{M}_{4\times 8}=\mathbf{U}_{4\times 4}\boldsymbol{\sigma}_{4\times 8}\mathbf{V}^T_{8\times 8}
\end{align}
where $\boldsymbol{\sigma}=diag(\sigma_1,..,\sigma_r)$ is diagonal matrix of singular values of $\mathbf{M}$ and the columns of $\mathbf{U}$ and $\mathbf{V}$ are left and right singular vectors of $\mathbf{M}$ respectively. Since the nullspace of $\mathbf{M}$ is spanned by right singular vectors corresponding to the vanishing singular values of $\mathbf{M}$, we express $\mathbf{c}$ as
\begin{align}
\label{sol}
\mathbf{c}_{8\times 1}=\begin{bmatrix}\mathbf{V}[:,r+1] & . & . & . & \mathbf{V}[:,8]\end{bmatrix}_{8\times (8-r)}\boldsymbol{\kappa}_{(8-r)\times 1}
\end{align}
where $\boldsymbol{\kappa}$ is the vector of coefficients and $r$ is rank of $\mathbf{M}$. The objective now is to determine $\boldsymbol{\kappa}$. 
First, using \eqref{expression}, we obtain an expression for the gradient $\nabla \mathcal{S}(\mathbf{x})$ of $\mathcal{S}(\mathbf{x})$ as
\begin{align}
\nabla \mathcal{S}(\mathbf{x})
\label{second}
&=\overbrace{\begin{bmatrix}
yz & y & 0 & z & 1 & 0 & 0 & 0\\
xz & x & z & 0 & 0 & 1 & 0 & 0\\
xy & 0 & y & x & 0 & 0 & 1 & 0
\end{bmatrix}}^{\mathbf{H}(x,y,z)_{3\times 8}}\begin{bmatrix}\mathbf{V}[:,r+1] & . & . & . & \mathbf{V}[:,8]\end{bmatrix}_{8\times (8-r)}\boldsymbol{\kappa}_{(8-r)\times 1}
\end{align}
Let $\hat{\mathcal{S}}(\hat{\mathbf{x}})$ be corresponding definition on face $\hat{e}$ of reference element $\hat{E}$ of $\mathcal{S}(\mathbf{x})$ on face $e$ of actual element $E$. Then, from \eqref{hat2}, 
\begin{align}
\label{gradient1}
\nabla \mathcal{S}(\mathbf{x})
=(DF_E)^{-T}(\hat{\mathbf{x}})\,\hat{\nabla}\hat{\mathcal{S}}(\hat{e})
\end{align}  
where $\hat{\nabla}\hat{\mathcal{S}}(\hat{e})$ can be either $\begin{bmatrix}
1 & 0 & 0
\end{bmatrix}^T$, $\begin{bmatrix}
0 & 1 & 0
\end{bmatrix}^T$ or $\begin{bmatrix}
0 & 0 & 1
\end{bmatrix}^T$ depending on whether $\hat{e}$ is normal to $\hat{x}$, $\hat{y}$ or $\hat{z}$ axis.
Equating \eqref{second} and \eqref{gradient1} for all four vertices of $e\in E$, we get the following system of equations for $\boldsymbol{\kappa}_{(8-r)\times 1}$
\begin{align}
\label{solfinal}
\begin{bmatrix}
\mathbf{H}(x_1,y_1,z_1)\\
\mathbf{H}(x_2,y_2,z_2)\\
\mathbf{H}(x_3,y_3,z_3)\\
\mathbf{H}(x_4,y_4,z_4)
\end{bmatrix}_{12\times 8}\begin{bmatrix}\mathbf{V}[:,r+1] & . & . & . & \mathbf{V}[:,8]\end{bmatrix}_{8\times (8-r)}\boldsymbol{\kappa}_{(8-r)\times 1}=\mathbf{B}_{12\times 1}
\end{align}
where $\mathbf{B}$ is obtained as
\begin{align*}
\mathbf{B}[(i-1)*3+1\rightarrow i*3,1]=(DF_E)^{-T}(\mathbf{\hat{v}_i})\,\hat{\nabla}\hat{\mathcal{S}}(\hat{e})
\end{align*}
where $\mathbf{\hat{v}_i}$, $i=1,2,3,4$ on $\hat{e}\in \hat{E}$ is the corresponding definition of $\mathbf{v}_i$, $i=1,2,3,4$ on $e\in E$.
The solution $\boldsymbol{\kappa}$ of \eqref{solfinal} is substituted into \eqref{sol} to obtain $\mathbf{c}$, which is then substituted into \eqref{expression} to obtain the polynomial expression of $\mathcal{S}(\mathbf{x})$.

\bibliographystyle{unsrt} 
\bibliography{diss}

\begin{thebibliography}{10}

\bibitem{dana-2018}
S.~Dana, B.~Ganis, and M.~F. Wheeler.
\newblock A multiscale fixed stress split iterative scheme for coupled flow and
  poromechanics in deep subsurface reservoirs.
\newblock {\em Journal of Computational Physics}, 352:1--22, 2018.

\bibitem{Carlo-1978}
C.~A. Felippa.
\newblock Iterative procedures for improving penalty function solutions of
  algebraic systems.
\newblock {\em International Journal for Numerical Methods in Engineering},
  12(5):821--836, 1978.

\bibitem{M-1978}
M.~J.~D. Powell.
\newblock Algorithms for nonlinear constraints that use lagrangian functions.
\newblock {\em Mathematical Programming}, 14(1):224--248, 1978.

\bibitem{J.O-1985}
J.~O. Hallquist, G.~L. Goudreau, and D.~J. Benson.
\newblock Sliding interfaces with contact-impact in large-scale lagrangian
  computations.
\newblock {\em Computer Methods in Applied Mechanics and Engineering},
  51:107--137, 1985.

\bibitem{Jua-1985}
J.~C. Simo, P.~Wriggers, and R.~L. Taylor.
\newblock A perturbed lagrangian formulation for the finite element solution of
  contact problems.
\newblock {\em Computer Methods in Applied Mechanics and Engineering},
  50(2):163--180, 1985.

\bibitem{P-1985}
P.~Wriggers and J.~C. Simo.
\newblock A note on tangent stiffness for fully nonlinear contact problems.
\newblock {\em International Journal for Numerical Methods in Biomedical
  Engineering}, 1(5):199--203, 1985.

\bibitem{H-1989}
H.~Parisch.
\newblock A consistent tangent stiffness matrix for three-dimensional
  non-linear contact analysis.
\newblock {\em International Journal for Numerical Methods in Engineering},
  28(8):1803--1812, 1989.

\bibitem{Panayioti-1992}
P.~Papadopoulos and R.~L. Taylor.
\newblock A mixed formulation for the finite element solution of contact
  problems.
\newblock {\em Computer Methods in Applied Mechanics and Engineering},
  94(3):373--389, 1992.

\bibitem{P-1993}
P.~Papadopoulos and R.~L. Taylor.
\newblock A simple algorithm for three-dimensional finite element analysis of
  contact problems.
\newblock {\em Computers and Structures}, 46(6):1107--1118, 1993.

\bibitem{T-2000}
T.~W. McDevitt and T.~A. Laursen.
\newblock A mortar-finite element formulation for frictional contact problems.
\newblock {\em International Journal for Numerical Methods in Engineering},
  48(10):1525--1547, 2000.

\bibitem{nag-2001}
N.~El-Abbasi and K.~J. Bathe.
\newblock Stability and patch test performance of contact discretizations and a
  new solution algorithm.
\newblock {\em Computers and Structures}, 79(16):1473--1486, 2001.

\bibitem{becker-2003}
R.~Becker, P.~Hansbo, and R.~Stenberg.
\newblock A finite element method for domain decomposition with non-matching
  grids.
\newblock {\em ESAIM Mathematical Modelling and Numerical Analysis},
  37(2):209--225, 2003.

\bibitem{puso-2003}
M.~A. Puso and T.~A. Laursen.
\newblock Mesh tying on curved interfaces in 3d.
\newblock {\em Engineering Computations}, 20(3):305--319, 2003.

\bibitem{Michae-2004}
M.~A. Puso and T.~A. Laursen.
\newblock A mortar segment-to-segment contact method for large deformation
  solid mechanics.
\newblock {\em Computer Methods in Applied Mechanics and Engineering},
  193(6-8):601--629, 2004.

\bibitem{wriggersbook}
P.~Wriggers.
\newblock {\em Computational Contact Mechanics}.
\newblock Springer, 2nd edition, 2006.

\bibitem{P-2008}
P.~Wriggers and G.~Zavarise.
\newblock A formulation for frictionless contact problems using a weak form
  introduced by nitsche.
\newblock {\em Computational Mechanics}, 41(3):407--420, 2008.

\bibitem{simo1992augmented}
J~C Simo and T~A Laursen.
\newblock An augmented lagrangian treatment of contact problems involving
  friction.
\newblock {\em Computers \& Structures}, 42(1):97--116, 1992.

\bibitem{glowinski1989augmented}
Roland Glowinski and Patrick Le~Tallec.
\newblock {\em Augmented Lagrangian and operator-splitting methods in nonlinear
  mechanics}.
\newblock SIAM, 1989.

\bibitem{adeli1994augmented}
Hojjat Adeli and Nai-Tsang Cheng.
\newblock Augmented lagrangian genetic algorithm for structural optimization.
\newblock {\em Journal of Aerospace Engineering}, 7(1):104--118, 1994.

\bibitem{conn1991globally}
Andrew~R Conn, Nicholas~IM Gould, and Philippe Toint.
\newblock A globally convergent augmented lagrangian algorithm for optimization
  with general constraints and simple bounds.
\newblock {\em SIAM Journal on Numerical Analysis}, 28(2):545--572, 1991.

\end{thebibliography}
\end{document}